\documentclass[a4paper,11pt]{article}      

\usepackage[dvips]{graphicx} 
\usepackage{graphics} 
\usepackage[latin1]{inputenc}
\usepackage{amssymb}         
\usepackage{amsmath}         
\usepackage{amsfonts}
\usepackage{a4wide}
\usepackage[usenames,dvipsnames]{xcolor}
\usepackage{curves}

\title{\LARGE \bf  About strong string stability of a vehicle chain with time-headway control}

\author{Arash Farnam and Alain Sarlette
\thanks{Arash Farnam is with the ID Lab, Department of Electronics and Information Systems (ELIS), Faculty of Engineering and Architecture, Ghent University; Technologiepark Zwijnaarde 914, 9052 Zwijnaarde(Ghent), Belgium {\tt\small arash.farnam@ugent.be.},  Alain Sarlette.~is with ELIS, Ghent University, Belgium; and with the QUANTIC lab, INRIA Paris; 2 rue Simone Iff, 75012 Paris, France. {\tt\small alain.sarlette@inria.fr}}%
\thanks{This paper presents research results of the Belgian Network DYSCO (Dynamical Systems, Control, and Optimization), funded by the Interuniversity Attraction Poles Programme, initiated by the Belgian State, Science Policy Office.}}

\begin{document}
\maketitle
\thispagestyle{empty}
\pagestyle{empty}


\begin{abstract}                
This paper deals with the problem of string stability in a chain of acceleration-controlled vehicles, i.e.~how input disturbances affect the distributed system for very long chains. There exist variants of string stability, like avoiding that a local disturbance gets amplified along the chain, or more strongly ensuring that the output vector's $p-$norm remains bounded for any bounded vector of input disturbances independently of the string length. They are all impossible to achieve with any linear controller if the vehicles only use relative information of few vehicles in front. Previous work has shown that adding absolute velocity into the controller, allows to at least avoid amplification of a local disturbance. In this paper, we consider the stronger definitions of string stability, under this same relaxation of using absolute velocity in the controller. We prove that the influence from input vector to output vector cannot be bounded independently of chain length in the most popular $2-$norm sense, with any bounded stabilizing linear controller; while a proportional derivative (PD) controller can guarantee it in the practically relevant $\infty-$norm sense. Moreover, we identify the disturbance acting on the leader as the main issue for string stability.
\end{abstract}


\section{Introduction}\label{sec:intro}

Platooning of vehicles is a method for increasing the capacity of roads and reducing fuel consumption, see e.g.~automated highway systems \cite{1}. Platoon coordination control must prevent the vehicles from too sudden acceleration/braking or even colliding, while maintaining short inter-vehicle distances for better performance. The most fundamental platoon is the vehicle chain, where all vehicles are aligned after each other. During the recent years, numerous works have considered different control strategies to stabilize each vehicle at a desired distance from its predecessor and follower in such chains \cite{1,3,4,Dirk}. The maybe unexpected challenge is the early observation that some bad behavior cannot be overcome with \emph{any} linear controller that feeds back relative distances between the vehicles.

More precisely, \emph{string instability} is a situation where the spacing error between consecutive vehicles grows unbounded when the number of vehicles increases to infinity, and string stability is the situation where this is avoided. (We intentionally keep the definition of ``spacing error'' loose at this point as there are several versions, to be detailed below.) This concept has spurred a lot of discussion and research since its definition in \cite{6,7}. Basically, it is known since \cite{6,7} that string stability cannot be achieved in a homogeneous chain of interconnected second-order integrators (e.g.~acceleration-controlled vehicles), with \emph{any} controller that is linear and whose local control actions are determined from the relative distance to a few directly preceding vehicles. This has attracted attention as a prototypical, unavoidable shortcoming of linear systems \cite{3}, \cite{4}. When each vehicle only reacts to its immediate predecessor, a straightforward proof of string instability follows from the Bode integral theorem \cite{25}. Indeed, the transfer function from error on vehicle $i-1$ to error on the following vehicle $i$, takes the form of a complementary sensitivity function, which unavoidably amplifies some frequencies of the disturbance \cite{7}. 

To investigate in more detail this problem and options to solve it, a distinction among several string stability notions was made. A weaker version on which researchers have concentrated first \cite{15,16,17,18,19,20}, and which we here label $L_2$ string stability (see Section \ref{sec:setting}), is to avoid that a \emph{single, local} $L_2$-bounded disturbance signal would have unbounded effects far off in the chain. The Bode integral argument indeed shows that particular local disturbances \emph{would} grow unbounded, unavoidably, with the simpler controllers. Solving this problem is a prerequisite for stronger versions, here labeled $(L_2,\ell_p)$ string stability, where one requests that a bounded vector disturbance signals, as the chain becomes infinitely long, induces a bounded vector of inter-vehicle distance errors in the same norm. For instance, $(L_2,\ell_\infty)$ string stability would request that if \emph{all} the vehicles are subject to an $L_2$-bounded input disturbance, then each inter-vehicle distance error should remain $L_2$-bounded. This appears to be the most ``practical'' formulation, at least when sticking to the traditional $L_2$ norm for signals in time. An $(L_\infty,\ell_\infty)$ version, corresponding to BIBO string stability, would arguably be even more significant in practice; but this would depart too much from the existing literature, and we leave it for future research in order to provide a significant novel point with $L_2$ signal norm. The $(L_2,\ell_p)$ string stability is the version considered originally by \cite{7}, just with input disturbances replaced by initial conditions. In standard work \cite{9,11,MidRev}, the $(L_2,\ell_2)$ version has in fact been the most popular proxy, for the benefit of easier analysis; it is called ``general $L_2$'' string stability in the review paper \cite{MidRev}.

Several lines of work have shown that $L_2$ string stability can be solved by adding a sufficiently strong feedback term proportional to \emph{absolute} velocity -- thus slightly enlarging the setting compared to purely relative information. The absolute velocity feedback can be obtained from a natural drag force \cite{Dirk}, although this would be less in line with developing ever more fuel efficient transportation means. In a somewhat subtler way, the absolute velocity term can also be obtained from a so-called \emph{time headway policy}, where the desired distance from a vehicle to its predecessor increases with the vehicle's velocity \cite{15,16}. This has the advantage to not favor slow vehicles, but it makes the effective distance velocity-dependent, where absolute velocity of the chain is a priori uncontrolled; so the interfacing with velocity control would have to be carefully investigated. It also remains to be seen exactly how (accurately) the absolute velocity is obtained in practice. While this absolute velocity solution has gathered serious attention as solving $L_2$ string stability \cite{15,16,17,18,19,20}, 
 it appears that no result so far has established its power for the stronger yet practically important $(L_2,\ell_p)$ versions. Those have only been investigated with even more information, e.g.~controllers relying on absolute position and/or on non-deteriorated knowledge, throughout the chain, of the leader's velocity profile \cite{26}. 

Therefore, the possibilities for time-headway to satisfy $(L_2,\ell_p)$ string stability have remained, somewhat surprisingly, open to date. Establishing these results is precisely the purpose of the present paper. We have both a positive result -- characterizing a PD controller which satisfies the ``practical'' $(L_2,\ell_\infty)$ string stability as requested in \cite{26}; and a negative result -- suggesting why these results were missing, namely because the more standard $(L_2,\ell_2)$ string stability notion cannot be satisfied by any controller that has bounded DC gain. We furthermore track the $(L_2,\ell_2)$ problem down to the effect of the leading vehicle only. This specificity of the string stability issue with time headway might come as a little surprising, and motivate further research towards circumventing string instability in $(L_k,\ell_p)$ sense: (i) the actual relevance of different norm choices; (ii) other contexts than vehicles, featuring possibly noise models with a ``perfectly clean'' leader; and (iii) the careful use of unbounded controllers, like a PID, in presence of other noise sources and contingencies. Note that under undirected coupling, controllers with integrator terms were proven to be unstable \cite[Thm.1.1]{9}. The concrete contribution of this paper is thus to establish a simple way to satisfy the strong and practical $(L_2,\ell_\infty)$ string stability; and to clarify that if one truly wants the $(L_2,\ell_2)$ version, then a more careful analysis using unbounded controllers will be necessary.

The paper is organized as follows. Section \ref{sec:setting} formalizes the problem setting and different definitions of string stability. In section \ref{FinRes3} we prove the impossibility result for $(L_2,\ell_2)$ string stability with bounded controllers, while in section \ref{FinRes4} we show how $(L_2,\ell_\infty)$ string satbility is achieved with the same controllers. In section \ref{FinRes5} an example illustrates our positive result.


\section{Problem Setting}\label{sec:setting}


\subsection{Model of Vehicle Chain}

Consider $N+1$ vehicles, whose position along the road at time $t$ we denote by $x(t) = (x_0(t),x_1(t),x_2(t),...,x_N(t)) \in \mathbb{R}^{N+1}$, with index $0$ denoting the leader. The focus of this work lies on the \emph{relative} position of consecutive vehicles, while their absolute value remains free. More precisely, we assume that the control objective is to stabilize the subspace $\{ x \in \mathbb{R}^{N+1} : x_{i}=x_{i-1}-r$ for $i=1,2,...,N$ \} for some given desired inter-vehicle distance $r>0$. Note that $r>0$ implies that vehicle $i$ is behind vehicle $i-1$ when they move with positive velocity. The configuration error vector thus writes $e(t) = (e_1(t),e_2(t),...,e_N(t))$ with
\begin{eqnarray}\label{error function}
e_i=x_{i-1}-x_i-r \   ,    \    i=1,2,...,N \, .
\end{eqnarray} 
The value of $x_0$ can then be independently guided as e.g.~a trajectory tracking command. The $N$ vehicles are modeled as isolated pure double-integrators with $u_i$ and $d_i$ as acceleration control input and disturbance input, respectively:
\begin{eqnarray}\label{system}
\ddot{x}_i(t)=u_i+d_i \   ,    \    i=0,1,2,...,N \, .
\end{eqnarray} 
To stabilize $e_1 = e_2 = ... = e_N = 0$, each vehicle $i$ adapts $u_i$ as a function of observed information about its neighboring vehicles. We here consider a \emph{unidirectional nearest-neighbor chain}, in line with much of the literature, where the feedback controller $u_i$ can depend on the relative information of one vehicle in front e.g.~their relative position  $x_{i-1}-x_{i}$ and relative velocity $\dot{x}_{i-1}-\dot{x}_{i}$. In addition, we allow $u_i$ to depend on the absolute velocity $\dot{x}_i$ of the corresponding vehicle like in \cite{15,16,18} and other papers, but unlike e.g.~\cite{17} we do not allow vehicle-to-vehicle communication. We suppose that the leading vehicle is a virtual one and assume $u_0=0$ as customary.





\subsection{Different Definitions of String Stability}

There are several variants of string stability in the literature, as mentioned in the introduction. In the following, we consider three of them: so called $L_2$ norm, $(L_2,\ell_2)$ norm, and $(L_2,\ell_\infty)$ norm string stability. The $L_2$ norm \cite{15,16,17,18,19,20} and $(L_2,\ell_2)$ norm \cite{9,11,MidRev} versions are most standard in the literature. The last one is closer to realistic concerns, but it seems that it was only recently examined in \cite{26}, while \cite{MidRev} also mentions versions with $L_2$ replaced by $L_p$ regarding the norm over time. The $L_2$ norm of a time-dependent scalar signal is denoted $\Vert x_i(\cdot) \Vert = \sqrt{\int_{-\infty}^{+\infty} \vert x_i(t) \vert^2 dt }$. For a time-dependent \emph{vector} e.g.~$x(t)$ , a lower index will indicate the discrete norm used on the \emph{vehicle index} dimension: the $(L_2, \ell_2)$ norm is
$$\Vert x(\cdot) \Vert_2 = \sqrt{ \sum_{i=0}^N  \int_{-\infty}^{+\infty} \vert x_i(t) \vert^2 dt }\, $$
and the $(L_2,\ell_\infty)$ norm is
$\;\; \Vert x(\cdot) \Vert_{\infty}  = \max_i \,(\Vert x_i(\cdot) \Vert)$.
\vspace{2mm}

\noindent \textbf{Definition ($(L_2,\ell_2)$ String Stability)}: \emph{The vehicle chain is $(L_2, \ell_2)$ string stable if, with the closed-loop dynamics, for every $\epsilon > 0$ there exists $\delta > 0$ such that: $\Vert d(\cdot) \Vert_2 < \delta$ implies $\Vert e(\cdot) \Vert_2 < \epsilon$, uniformly for all $N=1,2,...\;$.}\vspace{2mm}

\noindent \textbf{Definition ($L_2$ String Stability)}: \emph{The vehicle chain is $L_2$ string stable if, with the closed-loop dynamics, for every $\epsilon > 0$ there exists $\delta > 0$ such that: $\Vert d(\cdot) \Vert_2 < \delta$ implies $\Vert e_i(\cdot) \Vert < \epsilon$, uniformly for all $N=1,2,...\;$ and for all $i \in \{1,2,...,N\}$.}
\vspace{2mm}

\noindent \textbf{Definition ($(L_2,\ell_\infty)$ String Stability)}: \emph{The vehicle chain is $(L_2,\ell_\infty)$ string stable if, with the closed-loop dynamics, for every $\epsilon > 0$ there exists $\delta > 0$ such that: $\Vert d_j(\cdot) \Vert < \delta$ for all $j\in\{0, 1, ..., N\}$ (equivalently $\Vert d(\cdot) \Vert_\infty < \delta$) implies $\Vert e_i(\cdot) \Vert < \epsilon$ for all $i$ (equivalently $\Vert e(\cdot) \Vert_\infty < \epsilon$), uniformly for all $N=1,2,...\;$.}
\vspace{2mm}

\noindent \textbf{Remark 1 (admissible disturbances)}: These variants of string stability sometimes restrict the structure of the disturbance vector, e.g.~assuming $d_j=0$ for all $j>0$ to model disturbance on the leader only \cite{13}, or the opposite. Disturbances on the leader are indeed special, both practically since this is the ``active'' boundary of the chain, and for analysis since the controller on the leading vehicle is different; we will see that some results can indeed differ. A most notable point is that standard string stability studies assume disturbances on the inputs and/or initial conditions, but not on the measurements.
\vspace{2mm}

In a nutshell, the focus of string stability is that the configuration error must be bounded \emph{uniformly in $N$}. The weaker notion is $L_2$ string stability, as it bounds the sum of disturbance inputs but requests a bounded effect just independently for each $e_i$. $L_2$ string stability is a necessary condition but not sufficient to guarantee the stronger versions: $(L_2,\ell_2)$ string stability, where the sum-of-squares of the $e_i$ must be bounded too; and $(L_2,\ell_\infty)$ string stability, where the $e_i$ are considered individually but also the input disturbances $d_i$ need to be bounded only individually. We would argue that the latter is closest to a realistic physical situation -- a still more physical property might be the full BIBO version where also the $L_\infty$ norm is taken over time, but this would depart too much from the well-established literature to be our focus here. A priori, the $(L_2,\ell_2)$ and $(L_2,\ell_\infty)$ string stability are not in a definite relation with respect to each other. 

Several papers have considered the impossibility of $L_2$ string stability under relative information feedback, and proved how alternative settings using e.g.~absolute velocity feedback \cite{Dirk,15,17,18,19,20} do allow to achieve $L_2$ string stability with appropriate tuning. It has also been proved that it is possible to achieve $L_2$ norm string stability based on relative measurements only, using symmetric bidirectional controllers \cite{13} in which vehicles react not only to their predecessor but also symmetrically to their follower, \emph{but only provided disturbances are assumed to act exclusively on the leader}. However, regarding $(L_2,\ell_2)$ string stability, as well as the less-studied $(L_2,\ell_\infty)$ string stability, the situation is more negative. In the symmetric bidirectional control setting of \cite{13}, it has been proved that $(L_2,\ell_2)$ norm string stability cannot be achieved using any linear symmetric bidirectional controllers, see \cite{9,11}. For the controllers with absolute velocity feedback, satisfying the stronger versions of string stability has remained open so far, and we precisely set to answer this point.

\subsection{Considered control situations}

The first observations of string instability (see e.g.~\cite{7}) were made when each control input $u_i$ is reacting just to the relative distance $e_i$ with the vehicle in front, and trying to stabilize a constant inter-vehicle distance $r$. A relatively simple proof allows to check that string stability in $L_2$ sense, and thus a fortiori the stronger variants, is impossible in this case, with any linear controller that avoids pole cancellation. We recall the proof here as some related computations will be used below. A common assumption for the whole paper is that the controller should not be based on perfect pole cancellation, i.e.~$K(0) \neq 0$. Let $u_i(s) = K(s) e_i(s)$. From \eqref{error function},\eqref{system}, the closed-loop equation for the $e_i$ writes
$$e_i = T(s)\, e_{i-1} + \frac{1}{s^2+K(s)}\, (d_{i-1}-d_{i})\; ,$$
with
$$T(s)=\frac{K(s)}{s^2+K(s)} = \frac{R(s)}{1+R(s)}$$
where $R(s)=K(s)/s^2$. This takes the form of a complementary sensitivity function. Now assume for example that there is only disturbance on the leading vehicle, so
$$e_i= T(s)^{i-1}\frac{1}{s^2+K(s)}d_0 \; .$$
To guarantee $L_2$ string stability, with $N$ unboundedly large, it is then necessary in particular that $\vert T(j\omega) \vert \leq 1$ at all frequencies $\omega$. One concludes that this is impossible for a stable system, from the statement of Bode's Complementary Sensitivity integral \cite{25} which we recall below.
\vspace{2mm}

\noindent \textbf{Proposition 1:} \emph{Assume that the loop transfer function $R(s)$ of a system has (at least) a double pole at $s=0$. If the associated feedback system is stable, then the complementary sensitivity function $T(s)=\frac{R(s)}{1+R(s)}$ must satisfy: 
$$\int_{0}^{\infty} ln\mid T(j\omega)\mid d\omega/\omega^2 = \pi \sum_k\frac{1}{q^{(T)}_k} \ge 0 \, ,$$
where $\{ q^{(T)}_k \}$ are the zeros of $R(s)$ in the open right half plane. In particular, if $\vert T(j\omega) \vert < 1$ at some frequencies, then necessarily $\vert T(j\omega) \vert > 1$ at other frequencies.}
\vspace{2mm}

So, no linear controller of the type $u_i(s) = K(s) e_i(s)$ can achieve string stability.\\

For that reason, researchers have investigated controllers that depend not only on the relative distance between the vehicles $x_{i-1}-x_i$ but also on the \emph{absolute velocity} $\dot{x}_i$ of the vehicle itself. Regarding implementation, this can appear from velocity damping \cite{Dirk} or from a so-called time-headway policy \cite{15} where desired inter-vehicle distance $r$ would depend on $\dot{x}_i$. In the latter case, the form $e_i = x_{i-1}-x_i-h\dot{x}_i -r_0$ has been proposed with constant $r_0>0$ and time-headway parameter $h>0$, as illustrated on Fig.1. The linear controller $u_i(s) = K(s) e_i(s)$, in presence of time-headway spacing policy, implies the closed-loop equation
\begin{eqnarray}
\label{eq:controller2} e_i&=&\frac{K(s)}{s^2+(1+hs)K(s)}e_{i-1}\\ 
\nonumber & & + \frac{1}{s^2+(1+hs)K(s)}(d_{i-1}-(1+hs)d_i)\; ,
\end{eqnarray}
for $i=2,3,...,N$, and $e_1=\frac{1}{s^2+(1+hs)K(s)} (d_0-(1+hs)d_1)$. This controller was motivated by the following result \cite{15}.\vspace{2mm}

\color{black}

\noindent \textbf{Proposition 2:} \emph{The norm at $s=j\omega$ of transfer function $T(s)=\frac{K(s)}{s^2+(1+hs)K(s)}$ in \eqref{eq:controller2} is $< 1$ at all frequencies $\omega\neq0$, and its $H_\infty$ norm equals $T(0)=1$, if and only if the following condition holds:
\begin{equation}
\label{eq:headway2}
h>\underset{\omega}{\max} \sqrt{K_R(j\omega) \left(2-\omega^2 K_R(j\omega)\right)} + \omega K_J(j\omega)
\end{equation}
where $K_R(j\omega) = \frac{1}{2}(\frac{1}{K(j\omega)}+\frac{1}{K(j\omega)^*})$, $K_J(j\omega) = \frac{1}{2j}(\frac{1}{K(j\omega)}-\frac{1}{K(j\omega)^*})$, and the maximization runs over all $\omega$ for which the argument of the square root is positive.}\vspace{2mm}

\noindent \emph{Proof:} We just write $1/|T(j\omega)|^2 = |-\omega^2/K(j\omega) + (1 + h j \omega)|^2 > 1$ and we group real and imaginary parts to finally isolate $h$. \hfill $\square$\\

For particular controllers one can get easier criteria. E.g.~for a PD controller $K(s) = b s + a$, it is not hard to see that if $a>2b^2$ then the right hand side in \eqref{eq:headway2} is decreasing with $\omega$, and one gets the simple condition $h>\sqrt{2/a}$.

Note that the system with time headway is not subject to the Bode Integral, because we have $T(s) = \frac{R(s)}{1+R(s)}$ with $R(s) = K(s)/(s^2 + h s K(s))$ having a \emph{single} pole at $s=0$.




\section{Impossibility of $(L_2,\ell_2)$ string stability using bounded linear controllers with time-headway}\label{FinRes3}
\vspace{5mm}

We now consider the impossibility of achieving $(L_2,\ell_2)$ string stability using any \emph{bounded} stabilizing controller $K(s)$, in particular any controller satisfying $|K(0)| < \infty$, even in presence of time headway. (As we just recalled, without time headway i.e.~for $h=0$, it is already impossible to just achieve $L_2$ string stability.) We do this in two steps to identify that the main culprit is the disturbance on the leading vehicle: in essence, we can avoid that it gets amplified, but we cannot damp it fast enough along the chain with any bounded linear controller. In the literature, disturbances are mostly expected either everywhere, or exclusively on the leader like in \cite{13}. However, from an academic research point of view, it might be useful to know that only disturbances acting on the leader are causing the problem. The restriction to bounded controllers reflects the existing literature and avoids discussing other, unmodeled issues that might arise when e.g.~an integral term is present. In any case, the general conclusion may explain why a result about more than $L_2$ string stability was still missing regarding controllers with time-headway.



\subsection{No disturbance on leader, $d_0=0$}

While $d_0=0$ is not a practical situation, we treat it first to show, by linearity, that all problems essentially arise from $d_0$. We will show indeed that for $d_0=0$, one can achieve $(L_2,\ell_2)$ string stability using PD controllers with time headway.
\vspace{2mm}

\noindent \textbf{Theorem 1:}  \emph{There exists a pair $\,(K(s), \;h)\,$, where $h\geq 0$ is a sufficiently large constant time-headway satisfying Proposition 2 and $K(s)=bs+a$ is a stabilizing PD controller, such that the system \eqref{eq:controller2} is $(L_2,\ell_2)$ string stable provided $d_0=0$.}\vspace{2mm}

\emph{Proof:} The key point is to recognize that two effects of $d_i$ tend to compensate each other in $e_m$ with $m>i$. Indeed, we rewrite \eqref{eq:controller2} as
\begin{eqnarray*}
e_1&=&-L(s)d_1\\
e_i&=&-L(s)d_i+\sum_{m=2}^{i}T(s)^{i-m}P(s)d_{m-1} \; ,
\end{eqnarray*} 
with $P(s)=\frac{s^2}{(s^2+(1+hs)K(s))^2}$ and $L(s)=\frac{1+hs}{s^2+(1+hs)K(s)}$ and $T(s)$ defined as in Prop.2. We will assume that we place ourselves in the conditions of Prop.2, satisfying the related condition for a PD controller; it is not hard to check that the system is always stable with such PD controller. We next rewrite the dynamics in matrix form:
$$ e(s) = \big(-L(s) \mathbf{A}+ P(s) \mathbf{B}(s)\big) \, d(s)$$
with the $N \times (N+1)$ matrices
\begin{eqnarray} 
          \mathbf{A} &=& \begin{bmatrix}\
    0 & 1 & 0 & \dots  & 0 \\
    0 & 0 & 1& \dots  & 0 \\
        0 & 0 & 0& \dots  & 0 \\
    \vdots & \vdots & \vdots & \vdots & \vdots \\
     \\
    \nonumber 0 &0 & 0 & \dots &    1\\
\end{bmatrix} \;\; , \\[3mm]
    \mathbf{B}(s) &=&\begin{bmatrix}\
    0 & 0 & 0 & \dots  & 0 \\
    0 & 1 & 0& \dots  & 0 \\
        0 & T(s) & 1& \dots  & 0 \\
    \vdots & \vdots & \vdots & \vdots & \vdots \\
     \\
    \nonumber 0 &T(s)^{N-2} & T(s)^{N-3} & \dots &    0\\
\end{bmatrix} \;\; .
\end{eqnarray}    
We first use the triangle inequality to bound
\begin{eqnarray*}
\Vert e(s)\Vert_{2}&\le& \big( \vert L(s) \vert\, \Vert \mathbf{A}\Vert_2 + \vert P(s) \vert\, \Vert \mathbf{B}\Vert_2)  \Vert d(s) \Vert_2
\end{eqnarray*}     
with the induced matrix norms, i.e.~
$$\Vert \mathbf{D} \Vert_2 = \sqrt{\lambda_\text{max} (\mathbf{D}^{*} \mathbf{D})} \, $$
with $^*$ the complex conjugate transpose.
The proof now comes down to proving a bounded norm, independent of $N$ and $s=j\omega$, for the coefficient in front of $\Vert d \Vert_2$.

For the first term, since $\mathbf{A}^*\mathbf{A} = \text{diag}(0,1,1,1,...,1)$, we immediately have $\vert L(s) \vert\, \Vert \mathbf{A}\Vert_2 = \vert L(s) \vert$, and the latter can be bounded independently of $s=j\omega$ for a stable system.

For the second term, we obtain that the element $(m,n)$ of the matrix $\mathbf{B}^* \mathbf{B}$ equals
$$\;T(s)^{m-n}\sum_{j=0}^{N-m}\vert T(s) \vert^{2j}\;$$
for $m,n \in \{2,3,...,N\}$, $m\geq n$, symmetrically for $n>m$, and zero for the remaining terms.
The Gerschgorin circle theorem thus says that all the eigenvalues of $\mathbf{B}^* \mathbf{B}(j\omega)$ are comprised in the circles of respective center and radius
\begin{eqnarray*}
c^{(m)} &=& \sum_{j=0}^{N-m}\; |T(j\omega)|^{2j}\;\; ,\\
r^{(m)} &=& \left( \sum_{j=0}^{N-m}\; |T(j\omega)|^{2j} \right) \; \left( \sum_{n=2,n\neq m}^N \vert T \vert^{|m-n|} \right) \, .
\end{eqnarray*}
With a PD controller satisfying Proposition 2, we have $|T(j\omega)| < 1$ for all $\omega>0$ and we can bound each sum by the result of an infinite geometric series. This yields\vspace{2mm}

\noindent $\vert P(j\omega) \vert^2\, \Vert \mathbf{B}(j\omega)\Vert^2_2 \leq \vert P(j\omega) \vert^2 \max_m (c^{(m)}+r^{(m)}) $
\begin{eqnarray*}
&\leq& \frac{1}{1-|T(j\omega)|^2} \cdot \frac{2}{1-|T(j\omega)|} \cdot \vert P(j\omega) \vert^2 \\
&=& \frac{\vert L(j\omega)\vert^2}{1-\vert T(j\omega)\vert^2} \cdot \frac{2 \vert R(j\omega)\vert^2}{1-\vert T(j\omega)\vert}
\end{eqnarray*}
where $R(s) = \frac{s^2\,/\,(1+hs)}{s^2+(1+hs)K(s)}$. Every factor in this expression is bounded at large frequencies, so there just remains to investigate the limit at $\omega=0$. For the first factor we have
\begin{eqnarray*}
\frac{\vert L(j\omega)\vert^2}{1-\vert T(j\omega)\vert^2} &=& \frac{1}{|\text{-}\omega^2+(1+jh\omega)K(j\omega)|^2-|K(j\omega)|^2} \\
&&\simeq \frac{1}{a^2} \cdot \frac{1}{h^2\omega^2}
\end{eqnarray*}
for $\omega$ close to zero. For the second factor, we have $\vert R(s) \vert \simeq \frac{\omega^2}{a}$, while $\frac{1}{1-\vert T(j\omega)\vert}$ has a leading term of order $1/\omega$ at low frequencies. Thus in fact $\vert P(j\omega) \vert^2\, \Vert \mathbf{B}(j\omega)\Vert^2_2$ is of order $\omega^4 / \omega^3$ and converges to zero for low $\omega$. This gives a uniform bound on $\vert P(j\omega) \vert^2\, \Vert \mathbf{B}(j\omega)\Vert^2_2$ at all frequencies and thus concludes the proof.  \hfill $\square$

\subsection{Disturbance concentrated on $d_0$}

The chain's reaction to disturbances on the leader is slightly different, and we now show that this precludes the achievement of $(L_2,\ell_2)$ string stability with any linear controller of bounded DC gain.\vspace{2mm}

\noindent \textbf{Theorem 2:}  \emph{There exists no pair $\,(K(s) ,\;h)\,$, with $h\geq 0$ a constant time-headway and $K(s)$ a stabilizing controller with $K(0)$ finite, which would guarantee $(L_2,\ell_2)$ norm string stability of system \eqref{eq:controller2} when $d_0 \neq 0$.}
\vspace{2mm}

\noindent \emph{Proof:} We consider only a disturbance input $d_0$ that affects the leading vehicle, which leads to
\begin{eqnarray}\label{errors}
e_1&=&\frac{1}{s^2+(1+hs)K(s)}d_0\\
\nonumber e_i&=&T(s)^{i-1}\frac{1}{s^2+(1+hs)K(s)}d_0 \;\;\;, \;\; i\in\{2,3,...,N\}
\end{eqnarray}
with $T(s)$ defined as in Prop.~2. Then
\begin{eqnarray*}
\sum_{i=1}^{N} \vert e_i(s) \vert^2 &=& \sum_{i=0}^{N-1} \mid T(s)\mid^{2i} \cdot \frac{\vert d_0(s)\vert_2^2}{\vert s^2+(1+hs)K(s) \vert^2} \; .
\end{eqnarray*}
Take some $\beta>0$ and define $\alpha>0$ such that $\vert s^2+(1+hs)K(s) \vert^2\vert_{s=j\omega}<\alpha$ for all $\omega \in (-\beta,\beta)$. Now select any $\epsilon \in (0,\beta)$, and take an input disturbance concentrated at low frequencies such that
$$\int_{-\epsilon}^\epsilon \vert d_0(j\omega) \vert^2 d\omega \geq \frac{1}{2} \int_{-\infty}^{+\infty} \vert d_0(j\omega) \vert^2 d\omega \, .$$
Then
\begin{eqnarray}
\nonumber \Vert e(.) \Vert_2^2 \!\!&\!\! \geq \!\!&\!\! \int_{-\epsilon}^\epsilon \left. \sum_{i=0}^{N-1} \vert T(s)\vert^{2i} \frac{\vert d_0(s)\vert^2}{\vert s^2+(1+hs)K(s) \vert^2}\right\vert_{s=j\omega}  d\omega\\
\label{eq:2sumHA} \!\!&\!\! \geq \!\!&\!\! \frac{\Vert d_0(.) \Vert_2^2}{2} \, \frac{1}{\alpha} \;\;\; \min_{\omega \in (-\epsilon,\epsilon)} \sum_{i=0}^{N-1} \mid T(j\omega)\mid^{2i} \; .
\end{eqnarray}
Since $T(0)=1$, for any given $K(s)$ and $h$ and any $\delta>0$, there will always exist an $\epsilon$ such that
$\min_{\omega \in (-\epsilon,\epsilon)} \vert T(j\omega) \vert^2 > 1-\delta$. As $\delta$ can tend towards $0$ and $N$ towards infinity, the geometric sum in the second line of \eqref{eq:2sumHA} then cannot be bounded independently of $N$. \hfill $\square$\\

Theorem 2 implies the impossibility to achieve $(L_2,\ell_2)$ string stability in any cases where disturbances are expected \emph{at least} on the leading vehicle.As we mentioned earlier, this might explain why results in the literature are restricted to $L_2$ string stability, because the next-most popular setting would indeed be $(L_2,l_2)$ with bounded controllers $K(s)$. Luckily, there are two possible workarounds for this negative result. A first one is to allow $K(s)$ with unbounded DC gain, like a PID controller; we indeed show in \cite{ArashThesis} that a PID can be tuned to satisfy the $(L_2,l_2)$ definition of string stability with respect to input disturbances. However, the unbounded DC gain might require to investigate other effects more carefully, as unmodeled measurement noises or saturation effects could seriously deteriorate the situation. Another solution is to recognize that $(L_2, l_\infty)$ string stability might be a satisfactory achievement in practice. Indeed, for the latter case, we have the positive result that we present next.


\section{Satisfying $(L_2,\ell_\infty)$ String Stability with PD controller}\label{FinRes4}

We now turn to the positive part of the results, repeating a similar analysis to show how one does guarantee string stability in the ``practical'' $(L_2,\ell_\infty)$ sense using a PD controller with time headway.\vspace{2mm}

\noindent \textbf{Theorem 3:}  \emph{There exists a pair $\,(K(s), \;h)\,$, where $h\geq 0$ is a sufficiently large constant time-headway satisfying Proposition 2 and $K(s)=bs+a$ is a stabilizing PD controller, such that the system \eqref{eq:controller2} is $(L_2,\ell_\infty)$ string stable.}\vspace{2mm}

\noindent \emph{Proof:} The \emph{stability} of the system is easy to achieve with any positive tuning of the PD controllers as they lead to second-order polynomials in the denominators, also with time headway.

For \emph{string stability}, consider the worst case where there are disturbance inputs satisfying $\Vert d_i\Vert=\delta$ on all the vehicles $i\in\{0, 1, ..., N\}$. Following the notations of Thm.1 and Thm.2, we have
\begin{eqnarray}
\label{eq:nowgooderror} \Vert e_i\Vert &\le& \max_{s=j\omega} \left(  \sum_{m=2}^{i}|T^{i-m}P|+|T^{i-1}L_0|+| L| \right) \, \delta
\phantom{KkK}
\end{eqnarray}  
where $L_0(s)=\frac{1}{s^2+(1+hs)K(s)}$, while $T(s)=K(s)\,L_0(s)$, $L(s) = \frac{1+hs}{s^2+(1+hs)K(s)}$ and $P(s)=\frac{s^2}{(s^2+(1+hs)K)^2}$. By satisfying Proposition 2, we know that $\mid T(j\omega)\mid <1$ for all $\omega>0$, and then with the PD controller both $|L_0(s)| = |T(s)| / |a+bs|$ and $|L(s)| = |T(s)| \frac{|1+hs|}{|a+bs|}$ are bounded uniformly for all $s=j\omega$. The last two terms in \eqref{eq:nowgooderror} are thus bounded independently of $i$ and of $N$.

For the remaining term, we have
\begin{eqnarray}\label{practical}
&&\sum_{m=2}^{i}\vert T^{i-m}(j\omega)P(j\omega) \vert \le \frac{1}{1-|T(j\omega)|} \cdot |P(j\omega)| \phantom{KkK}\\
&&=\nonumber \frac{1}{\mid -\omega^2+(1+hj\omega) K(j\omega)\mid-\mid K(j\omega)\mid} \cdot\\
&&\nonumber \phantom{KkKkKkKkKkKk} \frac{\mid-\omega^2\mid}{\mid-\omega^2+(1+hj\omega)K(j\omega)\mid}.
\end{eqnarray}
We first check its behavior at low frequencies. By Taylor expansion we find
$\frac{1}{\mid-\omega^2+(1+j\omega h)K(j\omega)\mid-\mid K(j\omega)\mid}\simeq \frac{1}{\omega a}$ and $\mid\frac{-\omega^2}{-\omega^2+(1+hj\omega)K(j\omega)}\mid\simeq\frac{\omega^2}{a}$. For $\omega=0$ thus, \eqref{practical} converges to $0$. At low frequencies $\omega>0$, the deviation from $0$ in the right-hand side of \eqref{practical} is independent of $i$ and of $N$, and this provides a bound independent of $i$ and $N$ for the left-hand side. For any given controller satisfying Proposition 2, it is thus straightforward to identify some $\omega_0>0$ such that 
$\sum_{m=2}^{i}\vert T^{i-m}(j\omega)P(j\omega) \vert \; < \; 1/a$ for instance, for all $\omega \in (-\omega_0,\omega_0)$. There remains to prove that the same term remains bounded independently, of $i$ and $N$, for all $\omega>\omega_0$. With the proposed PD controller, for any $\omega_0>0$, there exists $\alpha<1$ such that $\mid T(j\omega)\mid\le\alpha$ for all $\omega>\omega_0$; this is checked for instance by ensuring a monotone decreasing Bode amplitude diagram, as we will draw in the example below. Then we have, for all $\omega>\omega_0$, a uniform bound on $\frac{1}{1-|T(j\omega)|}< \frac{1}{1-\alpha}$ and also on $|P(j\omega)|=|T(j\omega)|^2 \cdot \vert \omega / K(j\omega) \vert^2$. Together, all this provides a uniform bound on the first term of \eqref{eq:nowgooderror} and concludes the proof.\hfill $\square$

 
\section{Simulation} \label{FinRes5}

\begin{figure}[!ht]
\centering
\includegraphics[width=0.45\textwidth, trim=7mm 0mm 7mm 0mm, clip=true]{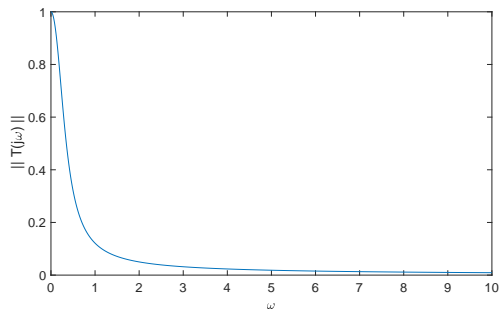}
\caption{Bode magnitude diagram of $T(s)$ for the example of Section \ref{FinRes5}. Note that we have avoided log-scales because it would squeeze the most important features, namely where $||T(j\omega)|| \simeq 1$.}\label{fig:2}
\end{figure}

Consider a PD controller $K(s) = bs + a$ for every vehicle with coefficients $a=b=1/6$. Taking $h=5$ satisfies the time-headway requirement given after Proposition 2. On Fig.\ref{fig:2} we show the Bode diagram of transfer function $T(s)$ and we can see that $T(s)$ is monotonically decreasing, as needed in the proof of Theorem 3.

\begin{figure}[!ht]
\centering
\includegraphics[width=0.5\textwidth, trim=6mm 0mm 7mm 0mm, clip=true]{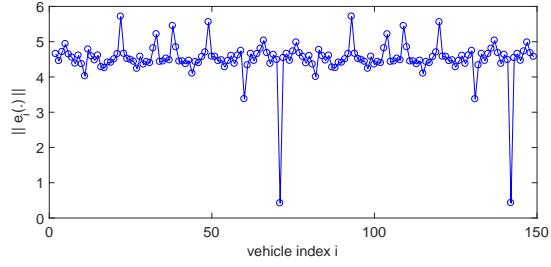}
\caption{$\Vert e_i \Vert$ as a function of $i$ when $L_2$-norm bounded disturbances are applied on all the vehicles. The $(L_2,\ell_\infty)$ string stability criterion for chain length $N$ is given by $\Vert e(.) \Vert_\infty = \max_{i\leq N}(\Vert e_i(.) \Vert)$; it remains bounded $\Vert e(.) \Vert_\infty < 6$ for any $N$.}\label{fig:3}
\end{figure}

\begin{figure}[!ht]
\centering
\includegraphics[width=0.5\textwidth, trim=7mm 0mm 7mm 1mm, clip=true]{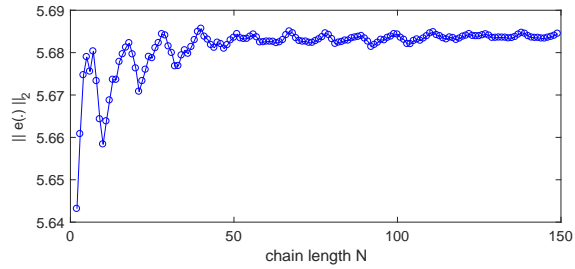}
\caption{Error $(L_2,\ell_2)$-norm $\Vert e(.)\Vert_2$ as a function of $N$ when the same disturbances are applied as on Fig.\ref{fig:3}, up to rescaling such that $\Vert d(.) \Vert_2 < 1$. This curve appears to saturate at large $N$.}\label{fig:4}
\end{figure}

To illustrate Thm.3 about $(L_2,\ell_\infty)$ string stability, we simulate this system for a string length $N=150$, applying pseudo-random $L_2$-norm bounded disturbances on all the vehicles; we just need to drop the last vehicles to get the behavior for shorter strings, thanks to the unidirectionality of the coupling. Figure \ref{fig:3} shows the $L_2$ norms of errors $e_i$ as a function of $i=1,2,...,150$. The error appears saturated, illustrating $(L_2,\ell_\infty)$ string stability, as expected from Thm.3. Repeating the simulation with other disturbance inputs, similarly normalized, we have always found a bound of about $\Vert e_i \Vert < 6$ for all $i$. Computing the $(L_2,\ell_2)$ string stability criterion with the same controller and same disturbance input, we must appropriately rescale $d$ as a function of $N$ to keep $\Vert d \Vert_2$ bounded. The corresponding criterion, shown on on Fig.\ref{fig:4}, appears to saturate and suggests string stability. However, this is due to the fact that our pseudo-random disturbance will give less and less weight to the leader disturbance $d_0$ as $N$ increases. Indeed, from Thms.1-2, we have seen that it is $d_0$ alone which causes the string instability, and this at low frequencies. 

\begin{figure}[!ht]
\centering
\includegraphics[width=0.5\textwidth, trim=7mm 0mm 7mm 1mm, clip=true]{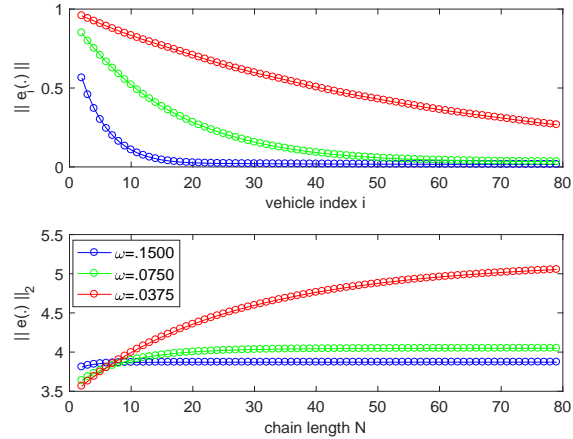}
\caption{Error $\Vert e_i \Vert$ as a function of $i$ (top) and $\Vert e(.)\Vert_2$ as a function of $N$ (bottom) for a disturbance input $d_0(t) = \sin(\omega t)$ and three values of $\omega$.}\label{fig:5}
\end{figure}

To illustrate the failure of $(L_2,\ell_2)$ string stability, as established by Thm.2, we must specifically apply a disturbance concentrated on the leader and at low frequencies. Fig.\ref{fig:5} shows the same string stability criteria as on Figs.\ref{fig:3} and \ref{fig:4}, but with the only disturbance $d_0 = \sin(\omega t)$ and three cases of different frequency $\omega$. As expected from the theorems, the $(L_2,\ell_\infty)$ criterion, always below 1, keeps decreasing along the chain since $\vert T(j\omega) \vert <1$ for all $\omega>0$. The $(L_2,\ell_2)$ criterion however grows unbounded as $N$ goes up and $\omega$ goes down, as the damping along the chain becomes weaker. Applying the same disturbance to another vehicle, e.g.~$d_2 = \sin(\omega t)$, we have observed that, as predicted by Thm.1, the $(L_2,\ell_2)$ criterion does remain bounded (figure not shown due to space constraints).

\section{Conclusion}

Controllers with absolute velocity feedback (e.g. time-headway) were known to solve the $L_2$ string stability issue for unidirectional vehicle chains. In this paper we have shown both theoretically and in simulations how two stronger definitions of string stability behave when using such controllers. From a more practical perspective, we have proved that a simple PD feedback with time-headway, does solve the $(L_2,\ell_\infty)$ version of string stability. Furthermore, we have pinned down a specific context for the impossibility to achieve $(L_2,\ell_2)$ string stability: the use of bounded controllers, as favored in the existing literature (e.g.~PD control); in combination with the dominating presence of low-frequency disturbance \emph{on the leader}, whose controller is subject to a boundary condition. This illustrates how subtle differences in definition can change the conclusions on string stability. Thus proxies on norms and error models are unfortunately not a great option and a careful study of the practically relevant criterion appears necessary for each application. It seems particularly relevant for future work to examine explicitly the BIBO-type version, namely $(L_\infty,\ell_\infty)$ string stability, which has gathered less attention in the literature than its $L_2$ proxies. Pushing this practical concern further, one could also seek to establish a precise quantitative tradeoff between absolute-velocity feedback $h$, control parameters in $K$, chain length $N$, and admissible local disturbances $\delta/\epsilon$, thus allowing to tune the controller to precision requirements $\delta$ for any finite $N$. This could enable to lower the value of $h$, gain some leeway towards attaining a better tradeoff with other control performance criteria, and possibly achieve stronger norms like $(L_2,\ell_2)$ bounds in this $N$-dependent-tradeoff sense.



\begin{thebibliography}{xx}

\bibitem{1}
K.~C.~Chu , ``Decentralized control of high-speed vehicular strings,'' \emph{Transportation Science}, Volume.~8, Pages.~361-384, 1974.

\bibitem{3} 
S.~Sheikholeslam and C.~Desoer, ``Longitudinal control of a platoon of vehicles,'' \emph{Proc. American Control Conf.},~Pages.~291-297, 1990.

\bibitem{4}
W.~Levine and M.~Athans, ``On the optimal error regulation of a string of moving vehicles,'' \emph{IEEE Trans. Automatic Control}, Volume.~11, Pages.~355-361, 1996.

\bibitem{Dirk}
J.A.~Rogge and D.~Aeyels, ``Vehicle Platoons Through Ring Coupling,'' \emph{IEEE Trans. Automatic Control}, Volume.~53, Pages.~1370-1377, 2008.

\bibitem{6}
D.~Swaroop and J.~Hedrick, ``String stability of interconnected systems,'' \emph{IEEE Trans. Automatic Control}, Volume.~41, Pages.~349-357, 1996.

\bibitem{7}
D.~Swaroop, ``String stability of interconnected systems: an application to platooning in automated highway systems,'' \emph{PhD thesis, University of California, Berkeley}, 1994.


\bibitem{9}
P.~Barooah and J.~P.~Hespanha, ``Error amplification and disturbance propagation in vehicle strings with decentralized linear control,'' \emph{Proc. IEEE Conf. on Decision and Control}, Pages~4964-4969, 2005.

\bibitem{11}
P.~Seiler, A.~Pant, and K.~Hedrick, ``Disturbance propagation in vehicle strings,'' \emph{IEEE Trans. Automatic Control}, Volume.~37, Pages. 1835 -1842, 1996.

\bibitem{13}
Y.~Yamamoto, and M.~C.~Smith, ``Bounded disturbance amplification for mass chains with passive interconnection,'' \emph{IEEE Trans. Automatic Control}, Volume.~47, Pages.~2534-2542, 2015.

\bibitem{15}   
S.~Klinge, and R.~H.~Middleton, ``Time headway requirements for string stability of homogenous linear unidirectionally connected systems,'' \emph{Proc. IEEE Conf. on Decision and Control}, Pages.~1992-1997, 2009.

\bibitem{16}
S.~Knorn, A.~Donaire, J.~C.~Aguero and R.~H.~Middleton, ``Passivity-based control for multi-vehicle systems subject to string constraints,"" \emph{Automatica}, Volume.~50, Pages.~3224-3230, 2014.

\bibitem{17}
 S.~Öncü, J.~Ploeg, N.~van de Wouw, and H.~Nijmeijer, ``Cooperative Adaptive Cruise Control:
network-aware analysis of string stability,'' \emph{IEEE Trans. Intelligent Transportation Systems}, Volume.~15, Pages.~1527-1537, 2014. 

\bibitem{18}
J.~Ploeg, N.~van de Wouw, and H.~Nijmeijer, ``$L_p$ string stability of cascaded systems: application to vehicle platooning,'' \emph{IEEE Trans. Control Systems Technology} Volume.~22, Pages.~786-793, 2014.

\bibitem{19}
J.~Ploeg, E.~S.~Kazerooni, G.~Lijster, N.~van de Wouw, and H.~Nijmeijer, ``Graceful degradation of Cooperative Adaptive Cruise Control,''  \emph{IEEE Trans. Intelligent Transportation Systems}, Volume.~16, Pages.~488-796, 2015. 

\bibitem{20}
V.~Milanés and S.~E.~Shladover, ``Modeling cooperative and autonomous adaptive cruise control dynamic responses using experimental data,'' \emph{Transportation Research Part C: Emerging Technologies}, V.~48, Pages.~285-300, 2014.


\bibitem{25}
D.~Looze and J.~Freudenberg, ``Tradeoffs and limitations in feedback systems,'' Ch.31 in \emph{The Control Handbook}, W. Levine, Ed. CRC Press, pp.~537-549, 1996.

\bibitem{26}
B.~Besselink and  K.~H.~Johansson, ``String Stability and a delay-based spacing policy for vehicle platoons subject to disturbances,'' \emph{IEEE Trans. Automatic Control}, Volume.~62, pp.~4376-4391, 2017.

\bibitem{MidRev}
S.~Stüdli, M.M. Seron and R.H. Middleton, ``From vehicular platoons to general networked systems: String stability and related concepts,'' \emph{Annual Reviews in Control}, 2017.

\bibitem{ArashThesis}
A.~Farnam, ``Towards impossibility and possibility results for string stability of platoon of vehicles,'' \emph{PhD Thesis, Ghent University,} 2018.

\end{thebibliography}
\end{document}